\documentclass[12pt]{amsart}
\usepackage[english]{babel}
\usepackage{amsmath}
\usepackage{amsthm}
\usepackage{amssymb}
\usepackage{mathrsfs}
\usepackage{enumerate}
\usepackage[notcite, final, notref]{showkeys}
\usepackage{amsfonts}
\usepackage{dsfont}
\usepackage{tikz}
\usepackage[siunitx]{circuitikz}
\usetikzlibrary{arrows,calc,chains,shapes,dsp}
\def\C{\mathbb C}
\def\R{{\mathbb R}}
\newtheorem{Pa}{Paper}[section]
\newtheorem{Tm}[Pa]{{\bf Theorem}}
\newtheorem{La}[Pa]{{\bf Lemma}}

\newtheorem{Rk}[Pa]{{\bf Remark}}
\newtheorem{Pn}[Pa]{{\bf Proposition}}

\newtheorem{Ex}[Pa]{{\bf Example}}
\newtheorem{Dn}[Pa]{{\bf Definition}}
\usepackage{xcolor}

\title[Passive LTI Systems: A Unified Framework]
{A Unified Framework for Continuous/Discrete\\[0.2cm]
Positive/Bounded Real State-Space Systems}
\author[I. Lewkowicz]{Izchak Lewkowicz}
\address{School of Electrical and Computer Engineering 
Ben-Gurion University of the Negev\\ P.O.B. 653\\ Beer-Sheva, 84105\\
Israel}
\email{izchak@bgu.ac.il}

\begin{document}
\bibliographystyle{plain}
\begin{abstract}
There are four variants of passive, linear time-invariant systems,
described by rational functions: Continuous or Discrete time, 
Positive or Bounded real. By introducing a quadratic matrix inequality
formulation, we present a unifying framework for state-space
characterization (a.k.a. Kalman-Yakubovich-Popov Lemma) of the above
four classes of passive systems.

\noindent
These four families are matrix-convex as rational functions, and
a slightly weaker version holds for the corresponding balanced,
state-space realization arrays.
\end{abstract}
\maketitle


\date{today}
\tableofcontents

\bibliographystyle{plain}
\section{Introduction}
\setcounter{equation}{0}

\noindent
In the study of dynamical systems, passivity is a fundamental property.
Thus, it has been extensively addressed in various frameworks. We here
focus on finite-dimensional, linear, time-invariant, passive systems
described by matrix-valued real rational functions of a complex
variable $z$. 
\vskip 0.2cm

\noindent
We shall use the following notation: $\C_L$ or $\C_R$ is the open Left
or Right half of the complex planes ($\overline{\C}_R$ is the closed
right half plane). Let also \mbox{$\mathbb D=\{z\in\C~:~ 1>|z|~\}$},
\mbox{$\overline{\mathbb D}=\{z\in\C~:~ 1\geq |z|~\}$}, be the 
open, closed unit disk and
\mbox{${\overline{\mathbb D}}^c=\{z\in\C~:~ |z|>1~\}$} is the 
the exterior of the closed unit disk\begin{footnote}{The
superscript $c$ stands for ``complement".}\end{footnote}.
\vskip 0.2cm

\noindent
For simplicity of exposition we begin with scalar functions terminology:
\begin{itemize}

\item[$(\alpha)$~]{} $\mathcal{P}$, Positive-Real (continuous-time)
analytically mapping $\C_R$ to its closure, $\overline{\C}_R$. See e.g.
\cite[Chapter 5]{AnderVongpa1973}, \cite[Chapter 7]{Belev1968},
\cite[Subsection 2.7.2]{BGFB1994}, \cite{DickDelsGenKam1985}. and 
\cite{MehrVanDo2020a}.

\item[$(\beta)$~]{} $\mathcal{B}$, Bounded-Real, (continuous-time) analytically
mapping, $\C_R$ to $\overline{\mathbb D}$. See \cite{AnderMoore1968},
\cite[Section 7.2]{AnderVongpa1973}, \cite[Chapter 7]{Belev1968},
\cite[Subsection 2.7.3]{BGFB1994}, and \cite{DickDelsGenKam1985}. 

\item[$(\gamma)$~]{} $\mathcal{DP}$, Discrete-Time-Positive-Real analytically
mapping ${\overline{\mathbb D}}^c$ to $\overline{\C}_R$, the closed right-half
plane. See e.g.  \cite{GoetVanGSuykVanDooDeMoor2003}, \cite{HitzAnder1969},
\cite{LiuXio2018},  \cite{MehrVanDo2020b} and \cite[Lemma 1]{XiaoHill1999}.

\item[$(\delta)$~]{} $\mathcal{DB}$, Discrete-Time-Bounded-Real analytically
mapping ${\overline{\mathbb D}}^c$ to $\overline{\mathbb D}$. See e.g.
\cite{Najs2013}, \cite{PremJury1994}, \cite{Vaidya1985} and 
\cite[Lemma 7]{XiaoHill1999}.
\end{itemize}

\noindent
The above families, are common in Engineering circles, for example
$\mathcal{P}$ and $\mathcal{DB}$ are used to describe passive systems in
continuous and discrete-time set-up, respectively.
\vskip 0.2cm

\noindent
For completeness, we point out that in mathematical analysis
community there are additional sets:
\vskip 0.2cm

\noindent
Herglotz or Carath\'{e}odory functions analytically map $\mathbb{D}$ to
$\overline{\C}_R$. See \cite{herglotz}. In other words, if $F(z)$ is a
Herglotz function, $F\left(z^{-1}\right)$ is a $\mathcal{DP}$ function. 
\vskip 0.2cm

\noindent
Schur functions analytically map $\mathbb{D}$ to its closure
$\overline{\mathbb D}$. See e.g. \cite{deBraRovn1966} and \cite{Schur}.
In other words, if $F(z)$ is a Schur function, $F\left(z^{-1}\right)$
is a $\mathcal{DB}$ function.
\vskip 0.2cm

\noindent
Recall that whenever $F(z)$ is an $p\times m$-valued rational function
with no pole at infinity, i.e. 
\mbox{$D:=\lim\limits_{z~\rightarrow~\infty}F(z)$}
is well-defined, one can associate with it a corresponding
\mbox{$(n+m)\times(n+m)$} state-space realization array, $R_F$ i.e.
\begin{equation}\label{eq:Realization}
F(z)=C(zI_n-A)^{-1}B+D\quad\quad\quad
R_F=\left({\footnotesize\begin{array}{c|c}A&B\\
\hline
C&D\end{array}}\right).
\end{equation}
The $(n+p)\times(n+m)$ realization $R_F$ in Eq. \eqref{eq:Realization}
is called {\em minimal}, if $n$ is the McMillan degree of $F(z)$.
\vskip 0.2cm

\noindent
In this work we focus on the case $p=m$ and examine characterizations of
the above four families of passive systems through the corresponding
state-space realizations. This is also known as the Kalman-Yakobovich-Popov
Lemma. For a (modest) account of the vast literature on the subject, 
beyond
those mentioned thus far, see e.g.
\cite{AlpayLew2011},
\cite{BarGohKaasRan2010}
\cite{CohenLew2007},
\cite{GusLik2006} (a survey),
\cite{Lewk2020a},
\cite{Lewk2020c},
\cite{Rant1996},
\cite{Will1972b} and
\cite{Will1976}.
For infinite-dimensional versions (all study Schur functions
in the above terminology) see e.g.
\cite{BallGroeterH2018},
\cite{BallStaffans2006},
\cite{Staffans2006}.
\vskip 0.2cm

\noindent
This work is organized as follows. In Section 2 we give the basic background.
Matrix-convex sets are introduced the Section 3. In Section 4, the main
result given and illustrated by a detailed example. In Section 5 we
specialize the main result to subsets of the families of functions discussed
in Section 4. In Section 6 the uniform framework is applied to show
matrix-convexity of systems, in the set-up of state-space realizations.

\section{Preliminary Background}
\label{Sec:Preliminaries}
\setcounter{equation}{0}

\noindent
In the sequel we shall denote by (${\mathbf H}_n$) $\overline{\mathbf H}_n$ 
the set of $n\times n$ (non-singular) Hermitian matrices. Skew-Hermitian matrices
are denoted by, $i\overline{\mathbf H}_n~$. It is common to take
$\overline{\mathbf H}$ and $i\overline{\mathbf H}$ as the matricial extension
of $\R$ and $i\R$, respectively. Then within Hermitian matrices
(${\mathbf P}_n$) $\overline{\mathbf P}_n$ will be the respective subsets of
positive (definite) semi-definite matrices. Recall that $\overline{\mathbf P}_n$
may be viewed as the closure of the open set
${\mathbf P}_n$.  
\vskip 0.2cm

\noindent
For ${\scriptstyle H}\in\mathbf{H}_n$ let us define the following sets
satisfying the Lyapunov inclusion,
\begin{equation}\label{eq:Lyapunov}
\begin{matrix}
\overline{\mathbf L}_H:=\left\{A\in\C^{n\times n} : 
\left(\begin{smallmatrix}A\\~\\I_n\end{smallmatrix}\right)^*
\left(\begin{smallmatrix}0&&H\\~\\H&&0\end{smallmatrix}\right)
\left(\begin{smallmatrix}A\\~\\I_n\end{smallmatrix}\right)
\in\overline{\mathbf P}_n\right\}
\\~\\
\mathbf{L}_H:=\left\{A\in\C^{n\times n} : 
\left(\begin{smallmatrix}A\\~\\I_n\end{smallmatrix}\right)^*
\left(\begin{smallmatrix}0&&H\\~\\H&&0\end{smallmatrix}\right)
\left(\begin{smallmatrix}A\\~\\I_n\end{smallmatrix}\right)
\in\mathbf{P}_n\right\},
\end{matrix}
\end{equation}
and the Stein inclusion,
\begin{equation}\label{eq:Stein}
\begin{matrix}
\overline{\mathbf S}_H:=\left\{A\in\C^{n\times n} : 
\left(\begin{smallmatrix}A\\~\\I_n\end{smallmatrix}\right)^*
\left(\begin{smallmatrix}-H&&0\\~\\~0&&H\end{smallmatrix}\right)
\left(\begin{smallmatrix}A\\~\\I_n\end{smallmatrix}\right)
\in\overline{\mathbf P}_n\right\}
\\~\\
\mathbf{S}_H:=\left\{A\in\C^{n\times n} : 
\left(\begin{smallmatrix}A\\~\\I_n\end{smallmatrix}\right)^*
\left(\begin{smallmatrix}-H&&0\\~\\~0&&H\end{smallmatrix}\right)
\left(\begin{smallmatrix}A\\~\\I_n\end{smallmatrix}\right)
\in\mathbf{P}_n\right\}.
\end{matrix}
\end{equation}
The above Quadratic Matrix Inclusion formulation is not the common way to
describe  the families $\mathbf{L}_H$ and $\mathbf{S}_H$. Yet it enables
us to present these sets in a common framework. This approach will be
taken a step forward in Theorem \ref{Tm:NewMain} below.
\vskip 0.2cm

\noindent
The sets $\overline{\mathbf L}_H$ and $\overline{\mathbf S}_H$ may be
viewed as the closure of the open sets $\mathbf{L}_H$ and $\mathbf{S}_H$,
respectively. The sets $\overline{\mathbf L}_H$ and $\mathbf{L}_H$ were
introduced and studied in \cite{CohenLew1997a}. In \cite{Ando2004}
Tsuyoshi Ando {\em characterized} the set $\mathbf{S}_H$.
\vskip 0.2cm

\noindent
We now resort to the classical Cayley transform. Recall that
$\mathcal{C}(A)$, the Cayley transform of a matrix $A\in\C^{n\times n}$,
is given by
\begin{equation}\label{eq:DefCayley}
\mathcal{C}(A):=(I_n-A)(I_n+A)^{-1}=-I_n+2(I_n+A)^{-1},
\quad\quad\quad-1\not\in{\rm spect}(A).
\end{equation}
Recall also that this transform is involutive, i.e.
whenever defined, $~\mathcal{C}\left(\mathcal{C}(A)\right)=A$.
\vskip 0.2cm

\noindent
It is well known that for a given ${\scriptstyle H}\in\mathbf{H}_n$,
\begin{equation}\label{eq:CayleyL}
\mathcal{C}\left(\overline{\mathbf L}_H\right)=
\overline{\mathbf S}_H\quad\quad\quad
\mathcal{C}\left(\mathbf{L}_H\right)=\mathbf{S}_H.
\end{equation}
When ${\scriptstyle-H}\in\mathbf{P}_n$, the set $\mathbf{L}_H$ is
associated with Hurwitz stability of differential equations of the
form \mbox{$\dot{x}=Ax$}. Similarly, when ${\scriptstyle H}\in\mathbf{P}_n$,
the set $\mathbf{S}_H$ is associated with Schur stability of difference
equations of the form \mbox{$x(k+1)=Ax(k)$}.
\vskip 0.2cm

\noindent
In the sequel, we shall focus on the special case where in Eqs.
\eqref{eq:Lyapunov} and \eqref{eq:Stein} $H=I_n$, i.e.
\begin{equation}\label{eq:L_I}
\begin{matrix}
\overline{\mathbf L}_{I_n}:=\left\{A\in\C^{n\times n}~:~
A+A^*\in\overline{\mathbf P}_n~\right\}
&&&
\overline{\mathbf S}_{I_n}:=\left\{A\in\C^{n\times n}~:~
1\geq\|A\|_2~\right\}.
\end{matrix}
\end{equation}
There is structural gap between $\mathbf{L}_H$ with $H\in\mathbf{H}_n$,
and the special case $\mathbf{L}_{I_n}$. This was addressed in
\cite[Theorems 2.4, 4.1]{Lewk2020a}. Similarly the structural gap
between $\mathbf{S}_H$, and the special case $\mathbf{S}_{I_n}$
was addressed in \cite[Theorem 2.1, Proposition 3.4]{Lewk2020c}.
\vskip 0.2cm

\noindent
We find it convenient to employ this matrix notation in order
to extend to matrix-valued set-up, the description of the above
four families of scalar real rational functions 

\begin{Dn}\label{Dn:Functions}
{\rm
Consider the following four families of $m\times m$-valued
real rational function.
\begin{itemize}

\item[$(\alpha)$~]{} $F\in\mathcal{P}$, means that
$\forall z\in\C_R$ one has that
$F(z)\in\overline{\mathbf L}_{I_m}$.
\vskip 0.2cm

\item[$({\beta})$~]{} $F\in\mathcal{B}$, means that
$\forall z\in\C_R$ one has that
$F(z)\in\overline{\mathbf S}_{I_m}$.
\vskip 0.2cm

\item[$(\gamma)$~]{} $F\in\mathcal{DP}$ means that
$\forall z\in\mathbb{D}$ one has that
$F(z)\in\overline{\mathbf L}_{I_m}$.
\vskip 0.2cm

\item[$(\delta)$~]{} $F\in\mathcal{DB}$ means that
$\forall z\in\mathbb{D}$ one has that
$F(z)\in\overline{\mathbf S}_{I_m}$.
\end{itemize}
\qed
}
\end{Dn}

\noindent
For completeness we next combine Definition \ref{Dn:Functions} along with
the Cayley transform. To this end recall that the Cayley transform of
\mbox{$m\times m$-valued} rational function $F(z)$, is defined as,
\[
\mathcal{C}(F)=(I_m-F)(I_m+F)^{-1}\quad\quad\quad
{\rm det}\left(F(z)+I_m\right)\not\equiv 0~~~\forall z\in\C .
\]
The functions in Definition \ref{Dn:Functions} satisfy the following
relations,
\begin{equation}\label{eq:CayleyFunctions}
\begin{matrix}
\mathcal{B}=\mathcal{C}\left(\mathcal{P}\right)
&&~&&
\mathcal{DB}=\mathcal{C}\left(\mathcal{DP}\right)
\\~\\
F(z)\in\mathcal{P}
&&\Longleftrightarrow&&
F\left({\scriptstyle\frac{1+z}{1-z}}\right)\in\mathcal{DP}
\\~\\
F(z)\in\mathcal{B}
&&\Longleftrightarrow&&
F\left({\scriptstyle\frac{1+z}{1-z}}\right)\in\mathcal{DB}.
\end{matrix}
\end{equation}
(In the Mathematical analysis terminology the formulation is more symmetric,
e.g. $F(z)$ belongs to $\mathcal{P}$ is equivalent to having
\mbox{$F\left({\scriptstyle\mathcal{C}(z)}\right)$} a Herglotz. function).

\section{Matrix-convex sets}
\label{Sec:MCsets}
\setcounter{equation}{0}

We next resort to the notion of a {\em matrix-convex} set, see e.g.
\cite{EffrWink1997} and more recently, \cite{Evert2018},
\cite{EverHeltKlepMcCull2018}, 
\cite{PassShalSol2018}.

\begin{Dn}\label{Dn:MatrixConvex}
{\rm
A family $\mathbf{A}$, of square matrices (of various
dimensions) is said to be}~ 
matrix-convex of level $n$,
{\rm if for
all $\nu=1,~\ldots~,~n$:\\
For all natural $k$, 
\[
\sum\limits_{j=1}^k\Upsilon_j^*\Upsilon_j=I_{\nu}
\quad\quad
\begin{smallmatrix}
\forall{\Upsilon}_j\in\C^{{\eta}_j\times\nu}
\\~\\
\eta_j\in[1,~\nu],
\end{smallmatrix}
\]
having $A_1,~\ldots~,~A_k$ (of various dimensions \mbox{$1\times 1$} through
\mbox{$\nu\times\nu$}) within $\mathbf{A}$, implies that also
\[
\sum\limits_{j=1}^k\Upsilon_j^*A_j\Upsilon_j~,
\]
belongs to $\mathbf{A}$.
\vskip 0.2cm

\noindent
If the above holds for all $n$, we say that the set $\mathbf{A}$ is}
matrix-convex.
\qed
\end{Dn}
\vskip 0.2cm

\noindent
In the rest of the section we briefly explore the notion of
matrix-convexity.

\begin{La}\label{La:MatrixConvex}~
\cite{Lewk2020a}.
The following sets are matrix-convex:
\begin{itemize}

\item[(i)~~~]{}
$\overline{\mathbf H},\quad\quad\quad
i\overline{\mathbf H},\quad\quad\quad
\overline{\mathbf P},\quad\quad\quad
{\mathbf P}$
\vskip 0.2cm

\item[(ii)~~]{}
$\{~A~:~{\rm Bound}\geq\| A\|_2\}~$ for some
${\rm Bound}>0$.
\vskip 0.2cm

\item[(iii)~]{}
The (open) closed set ($\mathbf{L}_I$) $\overline{\mathbf L}_I$,
see Eq. \eqref{eq:L_I}.
\end{itemize}
\end{La}

\noindent
Recall that the matrix-convexity condition is quite restrictive, so there
are not-too-many, non-trivial sets with this property. For example, the
sets (i) Toeplitz matrices, (ii) $~\{~A~:~{\rm Bound}\geq\| A\|_1\}~$
for some ${\rm Bound}>0$, (iii) $\mathbf{L}_P$ with 
\mbox{${\scriptstyle{\alpha}I\not=P}\in\mathbf{P}$}, are convex, but
not matrix-convex. Furthermore, matrix-convexity implies both classical
convexity and being unitarily-invariant, but the combination of these two
properties still falls short of characterizing matrix-convexity.
Indeed the set of positively scalar matrices, i.e of the form
${\alpha}I$, $\alpha>0$ is unitarily invariant and convex. However,
it is not matrix-convex: $A_1=\left(\begin{smallmatrix}2&&0\\~\\0&&2
\end{smallmatrix}\right)$ and $A_2=\left(\begin{smallmatrix}
3&&0\\~\\0&&3\end{smallmatrix}\right)$ belong to the set
but not the following combination
where 
\[
\underbrace{\left(\begin{smallmatrix}1&&0\\~\\0&&0
\end{smallmatrix}\right)}_{{\Upsilon}_1^*}
\underbrace{\left(\begin{smallmatrix}2&&0\\~\\0&&2
\end{smallmatrix}\right)}_{A_1}
\underbrace{\left(\begin{smallmatrix}1&&0\\~\\0&&0
\end{smallmatrix}\right)}_{{\Upsilon}_1}
+
\underbrace{\left(\begin{smallmatrix}0&&0\\~\\0&&1
\end{smallmatrix}\right)}_{{\Upsilon}_2^*}
\underbrace{\left(\begin{smallmatrix}3&&0\\~\\0&&3
\end{smallmatrix}\right)}_{A_2}
\underbrace{\left(\begin{smallmatrix}0&&0\\~\\0&&1
\end{smallmatrix}\right)}_{{\Upsilon}_2}
=
\left(\begin{smallmatrix}2&&0\\~\\0&&3\end{smallmatrix}\right)
\quad\quad\quad{\scriptstyle
{{\Upsilon}_1}^*{\Upsilon}_1+{{\Upsilon}_2}^*{\Upsilon}_2=I_2}~.
\]
Nevertheless, the four families of passive
rational functions we focus on, do share this property.

\begin{Pn}\label{Pn:ConvexRationalFunctions}
Each of the rational functions sets $\mathcal{P}$, $\mathcal{B}$, $\mathcal{DP}$
and $\mathcal{DB}$ is matrix-convex.
\end{Pn}

\noindent
{\bf Proof :}~ Let $F(z)$ be in $\mathcal{P}$ or in $\mathcal{B}$ and let
$F(z_o)$ be the image of a point $z_o$ which lies in the domain of interest
($\C_R$ and ${\overline{\mathbb D}}^c$ for $\mathcal{P}$ and $\mathcal{B}$,
respectively) From items $(\alpha)$, $(\gamma)$ in Definition \ref{Dn:Functions}
it follows that as a matrix, $F(z_o)$ is in $\overline{\mathbf L}_I$, see Eq.
\eqref{eq:L_I}, which is matrix-convex by item (iii) of Lemma \ref{La:MatrixConvex}.
\vskip 0.2cm

\noindent
In a similar way, let $F(z)$ be in $\mathcal{DP}$ or in $\mathcal{DB}$, and let
$F(z_o)$ be the image of a point $z_o$ which lies in the domain of interest
($\C_R$ and ${\overline{\mathbb D}}^c$ for $\mathcal{DP}$ and $\mathcal{DB}$,
respectively). From items $(\beta)$, $(\delta)$ in Definition \ref{Dn:Functions}
it follows that as a matrix, $F(z_o)$ is in $\overline{\mathbf S}_I$, see Eq.
\eqref{eq:L_I}, which is matrix-convex by item (ii) of Lemma \ref{La:MatrixConvex}.
\qed
\vskip 0.2cm

\noindent
In Proposition \ref{Pn:ConvexRealizationS} below, we offer a statement
analogous to Proposition \ref{Pn:ConvexRationalFunctions}, but in the
framework of realization arrays.
\vskip 0.2cm

\noindent
We end this section by pointing out that one can go beyond Proposition
\ref{Pn:ConvexRationalFunctions}.

\begin{Tm}\label{Tm:MaximalMatrixFunctions}
(I) \cite{Lewk2020a}.~The family $\mathcal{P}$, of $m\times m$-valued positive
real rational functions, is a cone, closed under inversion and a maximal
matrix-convex family of functions which is analytic in $\C_R~$.
\vskip 0.2cm

\noindent
Conversely, a maximal matrix-convex cone of $m\times m$-valued rational
functions, analytic in $\C_R$, containing the zero degree function
\mbox{$F(z)\equiv I_m$}, is the set $\mathcal{P}$.
\vskip 0.2cm

\noindent
(II) \cite{Lewk2020c}.~A family of $m\times m$-valued real rational
functions $F(z)$ which for all $z\in{\overline{\mathbb D}}^c$ is:\\
Analytic, matrix-convex and a maximal set closed under multiplication
among its elements, is the set $\mathcal{DB}$.
\vskip 0.2cm

\noindent
The converse is true as well.
\end{Tm}

\section{Characterization through State-Space: A Unified Framework}
\label{Sec:Framework}
\setcounter{equation}{0}

We start by recalling in the classical Kalman-Yakubovich-Popov type
characterization, through state-space realization, of the families in
Definition \ref{Dn:Functions}.

\begin{Tm}\label{Tm:ClassicalMain}
Consider the four families of $m\times m$-valued rational functions in
Definition \ref{Dn:Functions}. For each $l=\alpha$, $\beta$, $\gamma$
and $\delta$, let $R_{F_l}$ be a corresponding
state-space realization,
\[
F_l(z)=C_l(zI_n-A_l)^{-1}B_l+D_l~.
\]
(I)~ Assume Tet there exist $n\times n$ positive definite matrices $P$
so that the resulting \mbox{$(n+m)\times(n+m)$} matrices $Q_l$ are
positive semi-definite.

\begin{itemize}
\item[$(\alpha)$~]{}For
\[
Q_{\alpha}:=\left(\begin{smallmatrix}-PA-A^*P&&
C^*-PB\\~\\C-B^*P&&D+D^*\end{smallmatrix}\right),
\]
the function $F_{\alpha}(z)$ is Positive-Real. See e.g.
\cite{AlpayLew2011},
\cite[Chapter 5]{AnderVongpa1973},
\cite[Chapter 7]{Belev1968}, \cite[Subsection 2.7.2]{BGFB1994}, 
\cite{DickDelsGenKam1985},
and \cite[Theorem 3]{Will1972b}.

\item[$(\beta)$~]{}For 
\[
Q_{\beta}:=\left(\begin{smallmatrix}-PA-A^*P&&-PB\\~\\-B^*P&&I_m
\end{smallmatrix}\right)-\left(\begin{smallmatrix}C^*\\~\\D^*
\end{smallmatrix}\right)\left(\begin{smallmatrix}C&&D
\end{smallmatrix}\right),
\]
the function $F_{\beta}(z)$ is Bounded-Real. See e.g.
\cite{AnderMoore1968}, \cite[Section 7.2]{AnderVongpa1973}, 
\cite[Chapter 7]{Belev1968},
\cite[Subsection 2.7.3]{BGFB1994}, and \cite{DickDelsGenKam1985}

\item[$(\gamma)$~]{}For
\[
Q_{\gamma}:=
\left(\begin{smallmatrix}P&&C^*\\~\\C&&D+D^*\end{smallmatrix}\right)
-\left(\begin{smallmatrix}A^*\\~\\B^*\end{smallmatrix}\right)
\left(\begin{smallmatrix}A&&B\end{smallmatrix}\right),
\]
the function $F_{\gamma}(z)$ is Discrete-Time-Positive-Real.
See e.g. \cite{GoetVanGSuykVanDooDeMoor2003},
\cite{HitzAnder1969}, and \cite[Lemma 1]{XiaoHill1999}.

\item[$(\delta)$~]{}For,
\[
Q_{\delta}::=\left(\begin{smallmatrix}P-A^*PA&&-A^*PB\\~\\-B^*PA&&I_m
\end{smallmatrix}\right)-\left(\begin{smallmatrix}C^*\\~\\D^*
\end{smallmatrix}\right)\left(\begin{smallmatrix}
C&&D\end{smallmatrix}\right),
\]
the function $F_{\delta}(z)$ is Discrete-Time-Bounded-Real.  See
e.g. \cite{Najs2013}, \cite{PremJury1994}, \cite{Vaidya1985} and
\cite[Lemma 7]{XiaoHill1999}.
\end{itemize}
\vskip 0.2cm

\noindent
(II)~ In each of the four above cases, $l=\alpha$, $\beta$,
$\gamma$, $\delta$, if the state-space realization 
\mbox{$
F_l(z)=C_l(zI_n-A_l)^{-1}B_l+D_l$},
is minimal, then the converse is true as well.
\end{Tm}

\begin{Rk}
{\rm For completeness we recall in three extension of Theorem
\ref{Tm:ClassicalMain},
which are beyond the scope of this work.
\vskip 0.2cm

\noindent
{\bf 1.}~ If in each $Q_l$ in Theorem \ref{Tm:ClassicalMain}, having
${\scriptstyle P}\in\mathbf{P}_n$ is relaxed to 
${\scriptstyle H}\in\mathbf{H}_n$, then Generalized-positivity
(boundedness, ...) is obtained. For more details see \cite{AlpayLew2011},
\cite{AnderMoore1968},\cite[Theorem 10.2]{BarGohKaasRan2010} and
\cite{DickDelsGenKam1985}.
\vskip 0.2cm

\noindent
{\bf 2.}~ If each $Q_l$ in Theorem \ref{Tm:ClassicalMain} is restricted
to belong to $\mathbf{P}_{n+m}$, then the stronger notion of
Hyper-positivity (Hyper-boundednes, ...) is obtained. For further
details see \cite{AlpayLew2021}.
\vskip 0.2cm

\noindent
{\bf 3.}~ In \cite{AlpayLew2021} we addressed quantitative subsets of
$\mathcal{B}$, i.e. functions where,
\[
{\scriptstyle\sqrt{\frac{\eta-1}{\eta+1}}}
\geq\sup\limits_{z\in\C_R}\| F(z)\|_2\quad\quad
\eta\in(1,~\infty].
\]
A state-space characterization of the above $F$ belongs to this family when
$Q_{\beta}$ is substituted by 
\mbox{$Q_{\beta}(\eta)=
\left(\begin{smallmatrix}-PA-A^*P&&-PB\\~\\-B^*P&&I_m\end{smallmatrix}\right)
-{\scriptstyle\frac{\eta+1}{\eta-1}}\left(\begin{smallmatrix}C^*\\~\\D^*
\end{smallmatrix}\right)\left(\begin{smallmatrix}C&&D\end{smallmatrix}\right)
$}, which is positive-semidefinite.
\vskip 0.2cm

\noindent
Note that $\mathcal{B}$ is recovered when $\eta~\rightarrow~\infty$.
}
\qed
\end{Rk}

\noindent
To establish a unified framework, let us construct four $2(n+m)\times{2(n+m)}$
matrices, starting from
\begin{equation}\label{eq:Wdelta}
W_{\delta}=
\left(\begin{smallmatrix}
-P&&~0&&0&&0\\
~0&&-I_m&&0&&0\\
~0&&~0&&P&&0\\
~0&&~0&&0&&I_m
\end{smallmatrix}\right).
\end{equation}
Next taking a real symmetric orthogonal matrix
$U=\left(\begin{smallmatrix}\frac{1}{\sqrt{2}}I_n&0&\frac{1}{\sqrt{2}}I_n&0
\\0&I_m&0&0\\ \frac{1}{\sqrt{2}}I_n&0&-\frac{1}{\sqrt{2}}I_n&0\\0&0&0&I_m
\end{smallmatrix}\right)$, yields
\begin{equation}\label{eq:Wbeta}
U^*W_{\delta}U=W_{\beta}=
\left(\begin{smallmatrix}
~0&&~0&&-P&&0\\
~0&&-I_m&&~0&&0\\
-P&&~0&&~0&&0\\
~0&&~0&&~0&&I_m
\end{smallmatrix}\right).
\end{equation}
Now taking the permutation matrix
$U=\left(\begin{smallmatrix}0&&-I_{n+m}\\I_{n+m}&&~0\end{smallmatrix}\right)$,
one can obtain,
\begin{equation}\label{eq:Walpha}
W_{\delta}U=W_{\alpha}=\left(\begin{smallmatrix}~0&&0&&-P&&0\\~0&&0&&~0&&I_m\\
-P&&0&&~0&&0\\~0&&I_m&&~0&&0\end{smallmatrix}\right)
\end{equation}
and
\begin{equation}\label{eq:Wgamma}
W_{\beta}U=W_{\gamma}=\left(\begin{smallmatrix}-P&&0&&0&&0\\~0&&0&&0&&I_m\\
~0&&0&&P&&0\\~0&&I_m&&0&&0\end{smallmatrix}\right).
\end{equation}

\begin{Rk}\label{Rk:W}
{\rm
By construction the three real symmetric matrices: $W_{\alpha}$, $W_{\beta}$ 
and $W_{\gamma}$ are orthogonally similar to the real, block-diagonal,
symmetric matrix $W_{\delta}$. In particular the spectrum of each $W$ matrix,
is that of $\begin{smallmatrix}\pm{P}\end{smallmatrix}$ together with
$\begin{smallmatrix}\pm{I}_m\end{smallmatrix}~$.
}
\qed
\end{Rk}
\vskip 0.2cm

\noindent
The verification of the following, amounts to straightforward computation.

\begin{La}
Let $W_{\alpha}$, $W_{\beta}$, $W_{\gamma}$ and $W_{\delta}$ be as in Eqs.
\eqref{eq:Walpha}, \eqref{eq:Wbeta}, \eqref{eq:Wgamma} and \eqref{eq:Wdelta},
respectively then,
\begin{equation}\label{eq:Prototype}
Q_l=\left(\begin{matrix}R_{F_l}\\~\\I_{n+m}\end{matrix}\right)^*
W_l\left(\begin{matrix}R_{F_l}\\~\\I_{n+m}\end{matrix}\right)
\quad\quad\quad
l=\alpha, \beta, \gamma, \delta,
\end{equation}
with $Q_l$ as in Theorem \ref{Tm:ClassicalMain}.
\end{La}
\vskip 0.2cm

\noindent
This Lemma enables us now to cast the four items of Theorem
\ref{Tm:ClassicalMain} in a uniform framework.

\begin{Tm}\label{Tm:NewMain}
Let $R_{F_l}$ be an $(n+m)\times(n+m)$ realization of
$m\times m$-valued rational function, $F_l(z)$, 
\[
F_l(z)=C_l(zI_n-A_l)^{-1}B_l+D_l\quad\quad\quad
R_{F_l}=\left({\footnotesize\begin{array}{c|c}A_l&B_l\\
\hline
C_l&D_l\end{array}}\right)\quad\quad
l={\alpha}, {\beta}, {\gamma}, {\delta}.
\]
(I) With $W_{\alpha}$, $W_{\beta}$, $W_{\gamma}$ and $W_{\delta}$ from Eqs. 
\eqref{eq:Walpha}, \eqref{eq:Wbeta}, \eqref{eq:Wgamma} and \eqref{eq:Wdelta},
respectively consider the relation,
\begin{equation}\label{eq:Prototype1}
\underbrace{
\left(\begin{matrix}R_{F_l}\\~\\I_{n+m}\end{matrix}\right)^*
W_l
\left(\begin{matrix}R_{F_l}\\~\\I_{n+m}\end{matrix}\right)}_{
Q_l~~{\rm in~Eq.~\eqref{eq:Prototype}}}\in\overline{\mathbf P}_{n+m}~.
\end{equation}
Then the following is true:
\begin{itemize}
\item[$(\alpha)$~]{}If the condition in Eq. \eqref{eq:Prototype1} is satisfied
for $W_{\alpha}$ then $F_{\alpha}(z)$ is a Positive-Real function.
\vskip 0.2cm

\item[$(\beta)$~]{}If the condition in Eq. \eqref{eq:Prototype1} is satisfied
for $W_{\beta}$ then $F_{\beta}(z)$ is a Bounded-Real function.
\vskip 0.2cm

\item[$(\gamma)$~]{}If the condition in Eq. \eqref{eq:Prototype1} is satisfied
for $W_{\gamma}$ then $F_{\gamma}(z)$ is a Discrete-Time-Positive-Real function.
\vskip 0.2cm

\item[$(\delta)$~]{}If the condition in Eq. \eqref{eq:Prototype1} is satisfied
for $W_{\delta}$ then $F_{\delta}(z)$ is a Discrete-Time-Bounded-Real function.
\end{itemize}
\vskip 0.2cm

\noindent
(II)~ In each of the four above cases, if the realization $R_F$ is
minimal, then the converse is true as well.
\end{Tm}
\vskip 0.2cm

\noindent
We next illustrate an application of the unified framework in
Theorem \ref{Tm:NewMain}.
\vskip 0.2cm

\begin{Ex}
{\rm
We here address two forms of inversion and then
examine the behavior of the families $\mathcal{P}$ and
$\mathcal{DB}$ with respect to these two inversions.
\vskip 0.2cm

\noindent
We start with the classical one: When $F(z)$ is \mbox{$m\times m$-valued}
rational function, so that \mbox{${\rm det}(F(z))\not\equiv 0$}, one can
define \mbox{$\left(F(z)\right)^{-1}$} so that
\mbox{$\left(F(z)\right)^{-1}F(z)=I_m$} almost everywhere in $\C$.
\vskip 0.2cm

\noindent
Now, if \mbox{$F\in\mathcal{P}$}, by Theorem
\ref{Tm:MaximalMatrixFunctions}(I) $F^{-1}$ belongs to 
\mbox{$\mathcal{P}$} as well. In contrast, if
\mbox{$F\in\mathcal{DB}$}, from Theorem
\ref{Tm:MaximalMatrixFunctions}(II) it follows that $F^{-1}$ is in fact
an \mbox{``anti"-$\mathcal{DB}$} function in the sense that
\mbox{$\inf\limits_{z\in\overline{
\mathbb D}^c}{\sigma}_m\left((F(z))^{-1}\right)>1$}, where
${\sigma}_m(X)$ is the smallest singular value of a matrix
$X\in\C^{m\times m}$.
\vskip 0.2cm

\noindent
As a simple illustration take $f(z)=\frac{1}{2(2z+1)}$ (small letters
for scalar functions) and \mbox{$\left(f(z)\right)^{-1}=2(2z+1)$}.
Then, on the one hand both $f$ and $f^{-1}$ are in $\mathcal{P}$. Now
the same $f$ is in $\mathcal{DB}$ but \mbox{$f^{-1}$} belongs to
\mbox{``anti"-$\mathcal{DB}$}
in the sense that $|\left(f(z)\right)^{-1}|>1$ for all $|z|>1$.
\vskip 0.2cm

\noindent
One can now define another form of inversion, yielding an
\mbox{$m\times m$-valued} rational function (typically) of McMillan
degree $n$: Let $F(z)$ and $R_F$ be as in Theorem \ref{Tm:NewMain}. 
Treating the realization array $R_F$ as an \mbox{$(n+m)\times(n+m)$}
matrix, whenever it is non-singular. Then \mbox{$R_G=R_F^{-1}$} so that
\mbox{$R_GR_F=I_{n+m}~$}.
\vskip 0.2cm

\noindent
Let us return to the above simple illustration with 
\mbox{$f(z)=\frac{1}{2(2z+1)}$} so that its minimal realization is
\mbox{$R_f={\scriptstyle\frac{1}{2}}\left({\footnotesize\begin{array}
{r|r}-1&1\\ \hline 1&0\end{array}}\right)$} and then \mbox{$R_g=
R_f^{-1}={\scriptstyle 2}\left({\footnotesize\begin{array}{c|c}0&1\\
\hline 1&1\end{array}}\right)$}, which in turn means that
\mbox{$g(z)=\frac{2(2+z)}{z}~$}.
\vskip 0.2cm

\noindent
Now both \mbox{$f(z)=\frac{1}{2(2z+1)}$} and
\mbox{$g(z)=\frac{2(2+z)}{z}$} are positive real. However, while
$f(z)$ is also a $\mathcal{DB}$ function, $g(z)$ can not be a
$\mathcal{DB}$ function (nor ``anti"-$\mathcal{DB}$) since it maps
\mbox{$\{\R\smallsetminus[-1,~1]\}$} to \mbox{$(-2, 6)$}. 
\vskip 0.2cm

\noindent
Next, in the general framework, multiplying Eqs. \eqref{eq:Prototype},
\eqref{eq:Prototype1} by $R_G$ (recall $=\left(R_F\right)^{-1})$ from the
left and by ${R_G}^*$ from the right, yields
\begin{equation}\label{eq:Example}
\begin{matrix}
\underbrace{\begin{smallmatrix}{R_G}^*Q_lR_G\end{smallmatrix}}_{\tilde{Q}_l}
=
\begin{smallmatrix}{R_G}^*\end{smallmatrix}
\left(\begin{smallmatrix}R_F\\~\\I_{n+m}\end{smallmatrix}\right)^*
\begin{smallmatrix}W_l\end{smallmatrix}
\left(\begin{smallmatrix}R_F\\~\\I_{n+m}\end{smallmatrix}\right)
\begin{smallmatrix}R_G\end{smallmatrix}
=
\left(\begin{smallmatrix}I_{n+m}\\~\\R_G\end{smallmatrix}\right)^*
\begin{smallmatrix}W_l\end{smallmatrix}
\left(\begin{smallmatrix}I_{n+m}\\~\\R_G\end{smallmatrix}\right)
\\~\\
=
\left(\begin{smallmatrix}R_G\\~\\I_{n+m}\end{smallmatrix}\right)^*
\underbrace{
\left(\begin{smallmatrix}0&&I_{n+m}\\~\\I_{n+m}&&0\end{smallmatrix}\right)}_{U}
\begin{smallmatrix}W_l\end{smallmatrix}
\underbrace{
\left(\begin{smallmatrix}0&&I_{n+m}\\~\\I_{n+m}&&0\end{smallmatrix}\right)}_{U}
\left(\begin{smallmatrix}R_G\\~\\I_{n+m}\end{smallmatrix}\right).
\end{matrix}
\end{equation}
From Eq. \eqref{eq:Walpha} it follows that with
$U=\left(\begin{smallmatrix}0&&I_{n+m}\\~\\I_{n+m}&&0\end{smallmatrix}\right)$
one has that $U^*W_{\alpha}U=W_{\alpha}$ while $U^*W_{\delta}U=-W_{\delta}$.
In items (a) and (b) below,
we examine the system interpretation of this technical observation.
\vskip 0.2cm

\noindent
(a)~ If in Eq. \eqref{eq:Example} $F(z)$ is positive real,
i.e. $l={\alpha}$, one has that
\[
\left(\begin{smallmatrix}R_G\\~\\I_{n+m}\end{smallmatrix}\right)^*
\underbrace{
\left(\begin{smallmatrix}0&&I_{n+m}\\~\\I_{n+m}&&0\end{smallmatrix}\right)
\begin{smallmatrix}W_{\alpha}\end{smallmatrix}
\left(\begin{smallmatrix}0&&I_{n+m}\\~\\I_{n+m}&&0\end{smallmatrix}\right)
}_{U^*W_{\alpha}U=W_{\alpha}}
\left(\begin{smallmatrix}R_G\\~\\I_{n+m}\end{smallmatrix}\right)
=\begin{smallmatrix}\tilde{Q}_{\alpha}\end{smallmatrix}\in\overline{
\mathbf P}_{n+m}~, 
\]
namely,
\mbox{$
\left(\begin{smallmatrix}R_G\\~\\I_{n+m}\end{smallmatrix}\right)^*
\begin{smallmatrix}W_{\alpha}\end{smallmatrix}
\left(\begin{smallmatrix}R_G\\~\\I_{n+m}\end{smallmatrix}\right)
\in\overline{\mathbf P}_{n+m}$}. 
One can conclude that also $G(z)$ is positive real.
\vskip 0.2cm

\noindent
(b)~ If in Eq. \eqref{eq:Example} $F\in\mathcal{DB}$, i.e. $l={\delta}$, 
one has that
\[
\left(\begin{smallmatrix}R_G\\~\\I_{n+m}\end{smallmatrix}\right)^*
\underbrace{
\left(\begin{smallmatrix}0&&I_{n+m}\\~\\I_{n+m}&&0\end{smallmatrix}\right)
\begin{smallmatrix}W_{\delta}\end{smallmatrix}
\left(\begin{smallmatrix}0&&I_{n+m}\\~\\I_{n+m}&&0\end{smallmatrix}\right)
}_{U^*W_{\delta}U=-W_{\delta}}
\left(\begin{smallmatrix}R_G\\~\\I_{n+m}\end{smallmatrix}\right)=
\begin{smallmatrix}\tilde{Q}_{\delta}\end{smallmatrix},
\]
namely,
\mbox{$
\left(\begin{smallmatrix}R_G\\~\\I_{n+m}\end{smallmatrix}\right)^*
\begin{smallmatrix}(-W_{\delta})\end{smallmatrix}
\left(\begin{smallmatrix}R_G\\~\\I_{n+m}\end{smallmatrix}\right)
\in\overline{\mathbf P}_{n+m}~$}.
Thus, one can say that $G(z)$ is ``anti"-$\mathcal{DB}$: More
precisely,
\mbox{$\inf\limits_{z\in\mathbb{D}^c}{\sigma}_m\left((F(z))^{-1}\right)>1$}.
}
\qed
\end{Ex}
\vskip 0.2cm

\noindent
Consider the four families $\mathcal{P}$, $\mathcal{B}$, $\mathcal{DP}$ and
$\mathcal{DB}$. As already mentioned, as {\em rational functions} they
are related through the Cayley transform, see Eq. \eqref{eq:CayleyFunctions}. 
Theorem \ref{Tm:NewMain} suggests an additional inter-relations: Through
the corresponding {\em state-space realizations}. This is pursued in
Section \ref{Sec:Realizations}

\section{Two extreme cases}
\label{Sec:ExtremeCases}
\setcounter{equation}{0}

In this section we review two disjoint {\em subsets} of each of
the four families $\mathcal{P}$, $\mathcal{B}$, $\mathcal{DP}$
and $\mathcal{DB}$, which satisfy a Kalman-Yakubovich-Popov type
characterization in the form of Theorem \ref{Tm:NewMain}. 
\vskip 0.2cm

\noindent
To ease the reading, we present only samples of the four cases.

\subsection{Lossless}

Recall that the subset of $\mathcal{P}$ rational functions $F(z)$,
which in addition satisfy,
\[
\left(
\left(F(z)\right)^*
+
F(z)
\right)_{|_{z\in{i}\R}}=0,
\]
is called {\em Lossless Positive} (a.k.a. Positive Real Odd or
Foster), denoted by $\mathcal{LP}$.
\vskip 0.2cm

\noindent
In electrical circuits framework, $\mathcal{LP}$ functions are the
driving point impedance of $L-C$ network, i.e. a circuit comprised
of inductors and capacitors only (with zero resistance), see e.g.
\cite[Theorem 2.7.4]{AnderVongpa1973},
\cite[Ch 8, item 36]{Belev1968}, \cite[Section 4.2]{CohenLew2007}
and \cite{Lewk2020a}.
\vskip 0.2cm

\noindent
Following the relations in Eq. \eqref{eq:CayleyFunctions} in a similar
way, one may focus on the subset of $m\times m$-valued $\mathcal{B}$
rational functions $F(z)$, which in addition satisfy,
\[
{\left(F(z)\right)^*F(z)}_{|_{z\in{i}\R}}=I_m~,
\]
and call it Lossless Bounded, denoted by $\mathcal{LB}$.
\vskip 0.2cm

\noindent
Now, items $(\alpha)$ and $(\beta)$ of Theorem \ref{Tm:NewMain} takes the 
following form.

\begin{Pn}\label{Pn:Lossess}
Let $R_F$ be an $(n+m)\times(n+m)$ realization of
$m\times m$-valued rational function, $F(z)$, 
\[
F(z)=C(zI_n-A)^{-1}B+D\quad\quad\quad R_F=\left({\footnotesize\begin{array}
{c|c}A&B\\ \hline C&D\end{array}}\right).
\]
If,
\[
\left(\begin{matrix}R_F\\~\\I_{n+m}\end{matrix}\right)^*
W_{\alpha}
\left(\begin{matrix}R_F\\~\\I_{n+m}\end{matrix}\right)=0_{n\times n}~,
\]
with $W_{\alpha}$ as in Eq. \eqref{eq:Walpha}, then $F(z)$ is a
$\mathcal{LP}$ function.
\vskip 0.2cm

\noindent
If,
\[
\left(\begin{matrix}R_F\\~\\I_{n+m}\end{matrix}\right)^*
W_{\beta}
\left(\begin{matrix}R_F\\~\\I_{n+m}\end{matrix}\right)=0_{n\times n}~,
\]
with $W_{\beta}$ as in Eq. \eqref{eq:Wbeta}, then $F(z)$ is a
$\mathcal{LB}$ function.
\end{Pn}

\noindent
For further details see e.g. \cite{AlpayGoh1}, \cite{AlpayGoh2}
and \cite[p. 221-222, p.312]{AnderVongpa1973}.

\subsection{Hyper-Bounded}

In \cite{AlpayLew2021} we studied Hyper Bounded real rational:\\
A $m\times m$-valued rational function is said to be $\eta$-Hyper Bounded,
denoted by $\mathcal{HB}_{\eta}$, for some $\eta\in(1,~\infty]$ if,
\[
{\scriptstyle\sqrt{\frac{\eta-1}{\eta+1}}}\geq\sup\limits_{z\in\C_R}\|F(z)\|_2~.
\]
Equivalently (in the framework of Definition \ref{Dn:Functions}),
$F\in\mathcal{HB}_{\eta}$ if
\[
F(z)\in{\scriptstyle\sqrt{\frac{\eta-1}{\eta+1}}
}\cdot\overline{\mathbf S}_{I_m}\quad\quad\forall z\in\C_R~.
\]
This definition introduces the following partial ordering,
\[
\infty>\eta>\hat{\eta}>1\quad\quad\Longrightarrow\quad\quad
\mathcal{HB}_{\hat{\eta}}\subset\mathcal{HB}_{\eta}
\subset\{\mathcal{B}\smallsetminus\mathcal{LB}\}.
\]
Here, $W_{\beta}$ in Eq. \eqref{eq:Wbeta} is refined to,
\begin{equation}\label{eq:WbetaEta}
W_{{\beta},\eta}=\left(\begin{smallmatrix}~0&&~0&&-P&&0\\
~0&&\frac{1+\eta}{1-\eta}I_m&&~0&&0\\-P&&~0&&~0&&0\\
~0&&~0&&~0&&~~I_m\end{smallmatrix}\right)\quad\quad\quad
\eta\in(1,~\infty].
\end{equation}
Namely, $W_{\beta}=\lim\limits_{\eta~\rightarrow~\infty}W_{{\beta},\eta}~$.
\vskip 0.2cm

\noindent
Now, item $(\beta)$ of Theorem \ref{Tm:NewMain} takes the 
following form.

\begin{Pn}\label{Pn:HyperBounded}
Let $R_F$ be an $(n+m)\times(n+m)$ realization of
$m\times m$-valued rational function, $F(z)$, 
\[
F(z)=C(zI_n-A)^{-1}B+D\quad\quad\quad R_F=\left({\footnotesize\begin{array}
{c|c}A&B\\ \hline C&D\end{array}}\right).
\]
If,
\[
\left(\begin{matrix}R_F\\~\\I_{n+m}\end{matrix}\right)^*
W_{{\beta},\eta}
\left(\begin{matrix}R_F\\~\\I_{n+m}\end{matrix}\right)
\in\overline{\mathbf P}_{n+m}~,
\]
with $W_{{\beta},\eta}$ as in Eq. \eqref{eq:WbetaEta},
then $F(z)$ is a $\mathcal{HB}_{\eta}$ function.
\end{Pn}

\noindent
For details, see Lemma 2.3 and Eq. (2.6) in \cite{AlpayLew2021}.
\vskip 0.2cm

\noindent
For completeness, we point out that as before, one can definite
analogue sets by mimicking the relations in Eq.
\eqref{eq:CayleyFunctions}. For example, a $m\times m$-valued
rational function is said to be $\eta$-Hyper Discrete Bounded,
denoted by $\mathcal{HDB}_{\eta}$, for some $\eta\in(1,~\infty]$
if
\[
{\scriptstyle\sqrt{\frac{\eta-1}{\eta+1}}}\geq\sup\limits_
{z\in\overline{\mathbb D}^c}\|F(z)\|_2~.
\]
We shall not further pursue this direction.

\section{Sets of Matrix-convex Realization arrays}
\label{Sec:Realizations}
\setcounter{equation}{0}

\noindent
In this section we address inter-relations within {\em families of
realization arrays}~ associated with rational functions. As a
preliminary step we recall in the classical notion of transformation
of coordinates: Substituting a given state-space realization $R_F$ by 
\mbox{$
\left(\begin{smallmatrix}T^{-1}&0\\0~~&I_m\end{smallmatrix}\right)
R_F
\left(\begin{smallmatrix}T&~0\\0&~~I_m\end{smallmatrix}\right)$},
for some non-singular $n\times n$ matrix $T$.

\begin{La}\label{La:Ws}
Consider the framework of Theorem \ref{Tm:NewMain} for some
$~l\in\{{\alpha}, {\beta}, {\gamma}, {\delta}\}$.\\
Up to a
change of coordinates, one can take in Eq. \eqref{eq:Prototype},
\begin{equation}\label{eq:BalancedPrototype}
\left(\begin{matrix}R_{F_l}\\~\\I_{n+m}\end{matrix}\right)^*
\hat{W}_l
\left(\begin{matrix}R_{F_l}\\~\\I_{n+m}\end{matrix}\right)
\in\overline{\mathbf P}_{n+m}~,
\end{equation}
where the $\hat{W}$'s are associated with balanced realization, i.e.
\begin{equation}\label{eq:BalancedW}
\begin{matrix}
\hat{W}_{\alpha}&=&\left(\begin{smallmatrix}~0&&0&&-I_n&&0\\~0&&0&&~0&&I_m\\
-I_n&&0&&~0&&0\\~0&&I_m&&~0&&0\end{smallmatrix}\right)
&&&&
\hat{W}_{\beta}&=&\left(\begin{smallmatrix}~0&&~0&&-I_n&&0\\~0&&-I_m&&~0&&0\\
-I_n&&~0&&~0&&0\\~0&&~0&&~0&&I_m\end{smallmatrix}\right)
\\~\\
\hat{W}_{\gamma}&=&\left(\begin{smallmatrix}-I_n&&0&&0&&0\\~0&&0&&0&&I_m\\
~0&&0&&I_n&&0\\~0&&I_m&&0&&0\end{smallmatrix}\right)
&&&&
\hat{W}_{\delta}&=&\left(\begin{smallmatrix}-I_n&&~0&&0&&0\\~0&&-I_m&&0&&0\\
~0&&~0&&I_n&&0\\~0&&~0&&0&&I_m\end{smallmatrix}\right).
\end{matrix}
\end{equation}
\end{La}

\noindent
A system whose realization satisfies Eq. \eqref{eq:BalancedPrototype} with
$\hat{W}_{\alpha}$ from Eq. \eqref{eq:BalancedW} is called
``internally passive", see \cite[Definition 3]{Will1976}
\vskip 0.2cm

We also need to resort to following.

\begin{Dn}\label{Dn:n,mMatrixConvex}
{\rm
For all $k$, let $v_j\in\C^{(n+m)\times(n+m)}$,
$j=1,~\ldots~,~k$ be block-diagonal so that
\begin{equation}\label{eq:n,mIsometry}
\sum\limits_{j=1}^k
\underbrace{
\left(\begin{smallmatrix}{\Upsilon}_{j,n}&&0\\~\\0&&{\Upsilon}_{j,m}
\end{smallmatrix}
\right)^*}_{{\Upsilon}_j^*}
\underbrace{
\left(\begin{smallmatrix}{\Upsilon}_{j,n}&&0\\~\\0&&{\Upsilon}_{j,m}
\end{smallmatrix}
\right)}_{{\Upsilon}_j}
=
\left(\begin{smallmatrix}
I_n&&0\\~\\0&&I_m
\end{smallmatrix}
\right).
\end{equation}
A set $\mathbf{R}$ of $(n+m)\times(n+m)$ matrices is said to be~
$n,m$-{\em matrix-convex}~ if having
${\scriptstyle R_1,~\ldots~,~ R_k}$ in $\mathbf{R}$,
implies that also
\begin{equation}\label{eq:DefRgBR}
R_F:=\sum\limits_{j=1}^k
\underbrace{
\left(\begin{smallmatrix}{\Upsilon}_{j,n}&&0\\~\\0&&{\Upsilon}_{j,m}
\end{smallmatrix}\right)^*}_{{\Upsilon}_j^*}
\underbrace{
\left(\begin{smallmatrix}A_j&&B_j\\~\\C_j&&D_j\end{smallmatrix}
\right)}_{R_{F_j}}
\underbrace{
\left(\begin{smallmatrix}{\Upsilon}_{j,n}&&0\\~\\0&&{\Upsilon}_{j,m}
\end{smallmatrix}\right)}_{{\Upsilon}_j},
\end{equation}
belongs to ${\mathbf R}$
for all natural $k$ and all ${\Upsilon}_j\in\C^{(n+m)\times(n+m)}$.
\qed
}
\end{Dn}
\vskip 0.2cm

\noindent
In \cite{Lewk2020a} it was pointed out that the notion of
\mbox{$n,m$-{\em matrix-convexity}} is intermediate between (the more strict)
{\em matrix-convexity}, and (weaker) classical convexity.
\vskip 0.2cm

\noindent
We now pose the following question: For a natural parameter $k$, let 
$F_1(s)~,~\ldots~,~F_k(s)$ be a family of \mbox{$m\times m$-valued} rational
functions all from the same family, admitting \mbox{$(n+m)\times(n+m)$}
realizations, i.e.
\begin{equation}\label{eq:RGj}
R_{F_j}=\left({\footnotesize\begin{array}{l|r}\hat{A}_j&\hat{B}_j\\
\hline\hat{C}_j&\hat{D}_j\end{array}}\right)
\quad\quad\quad j=1,~\ldots~,~k.
\end{equation}
Let $R_F$, a realization of an \mbox{$m\times m$-valued} rational
function $F(z)$, be as in Eq. \eqref{eq:DefRgBR}. We now address the
following problem:\\
Under what conditions having
$F_1(z),~\ldots~,~F_k(z)$ in Eq. \eqref{eq:RGj} all in
$\mathcal{P}$ (or $\mathcal{B}$ or $\mathcal{DP}$ or $\mathcal{DB}$)
implies that also the resulting $F(z)$ in Eq. \eqref{eq:DefRgBR}
belongs to the same set?
\vskip 0.2cm

\noindent
If such a property holds this suggests that out of a small number of
``extreme points" of balanced realizations of $\mathcal{P}$ (or
$\mathcal{B}$ or $\mathcal{DP}$ or $\mathcal{DB}$) rational functions,
one can construct a whole ``matrix-convex-hull" realizations of
functions within the same family.  As a sample application, this may
enable one to perform a simultaneous balanced truncation model order
reduction of a whole family of bounded real functions, in the
spirit of \cite[Section 5]{CohenLew1997b}.
\vskip 0.2cm

\noindent
Before addressing this question, a word of caution: For example, 
$R_1=\left({\footnotesize\begin{array}{c|c}A&B\\ \hline C&D
\end{array}}\right)$ and $R_2=\left({\footnotesize\begin{array}{r|r}
A&-B\\ \hline-C&D\end{array}}\right)$ are two realization of the same
rational function. Furthermore, $R_1$ is minimal (balanced) if and
only if $R_2$ is minimal (balanced). However, \mbox{$R_3=
{\scriptstyle\frac{1}{2}}(R_1+R_2)=\left({\footnotesize\begin{array}
{c|c}A&0\\ \hline 0&D\end{array}}\right)$} is only a non-minimal
realization of a zero degree rational function $F(s)\equiv D$.
\vskip 0.2cm

\noindent
In a similar way, even when the ``extreme points" realizations in Eq.
\eqref{eq:RGj} are all balanced, the resulting realization $R_F$ in
Eq. \eqref{eq:DefRgBR}, may be not minimal.
\vskip 0.2cm

\noindent
We now return to the above question, 

\begin{Pn}\label{Pn:ConvexRealizationS}
Consider the framework of Lemma \ref{La:Ws} where $~l\in\{{\alpha}$,
${\beta}, {\gamma}, {\delta}\}$ is prescribed.
For a natural parameter $k$, let  $F_{1,l}(z)~,~\ldots~,~F_{k,l}(z)$
be a family of
\mbox{$m\times m$-valued} rational functions, admitting $(n+m)\times(n+m)$
realizations as in Eq. \eqref{eq:RGj}, satisfying all Eq.
\eqref{eq:BalancedPrototype} i.e.
\begin{equation}\label{eq:SimultaneousPrototype}
\left(\begin{matrix}R_{F_{j,l}}\\~\\I_{n+m}\end{matrix}\right)^*
\hat{W}_l
\left(\begin{matrix}R_{F_{j,l}}\\~\\I_{n+m}\end{matrix}\right)
=Q_{j,l}\in\overline{\mathbf P}_{n+m}~~~~\begin{matrix}
l\in\{\alpha,~\beta,~\gamma,~\delta\}~~{\rm prescribed}\\~\\
j=1,~\ldots~,~k.
\end{matrix}
\end{equation}
Then, $R_F$ in Eq. \ref{eq:DefRgBR} satisfies the same relation, i.e.
each of sets $\mathcal{P}$, $\mathcal{B}$, $\mathcal{DP}$ and
$\mathcal{DB}$ is a realization-$m,n$-matrix-convex.
\end{Pn}

\noindent
{\bf Proof :}~ Assume that 
Eq. \eqref{eq:SimultaneousPrototype} holds for $l=\alpha$, i.e.
\[
\underbrace{\left(\begin{smallmatrix}A_j&&B_j\\~\\C_j&&D_j
\\~\\I_n&&0\\~\\0&&I_m\end{smallmatrix}\right)^*
\overbrace{
\left(\begin{smallmatrix}~0&&0&&-I_n&&0\\~\\~0&&0&&~0&&I_m\\~\\
-I_n&&0&&~0&&0\\~\\~0&&I_m&&~0&&0\end{smallmatrix}\right)
}^{\hat{W}_{\alpha}}
\left(\begin{smallmatrix}A_j&&B_j\\~\\C_j&&D_j\\~\\I_n&&0
\\~\\0&&I_m\end{smallmatrix}\right)}_{Q_j}\quad\quad
\begin{smallmatrix}j=1,~\ldots~,~k\\~\\
Q_j\in\overline{\mathbf P}_{n+m},
\end{smallmatrix}
\]
and consider matrix-convex combination of realizations
as in Eq. \eqref{eq:DefRgBR},
\[
\begin{matrix}
\left(\begin{smallmatrix}
\sum\limits_{j=1}^k{\Upsilon}_{j,n}^*A_j{\Upsilon}_{j,n}&&
\sum\limits_{j=1}^k{\Upsilon}_{j,n}^*B_j{\Upsilon}_{j,m}\\~\\
\sum\limits_{j=1}^k{\Upsilon}_{j,m}^*C_j {\Upsilon}_{j,n}
&&
\sum\limits_{j=1}^k{\Upsilon}_{j,m}^*
D_j
{\Upsilon}_{j,m}
\\~\\
I_n&&0
\\~\\
0&&I_m\end{smallmatrix}\right)^*
\underbrace{
\left(\begin{smallmatrix}~0&&0&&-I_n&&0\\~\\~0&&0&&~0&&I_m\\~\\
-I_n&&0&&~0&&0\\~\\~0&&I_m&&~0&&0\end{smallmatrix}\right)
}_{\hat{W}_{\alpha}}
\left(\begin{smallmatrix}
\sum\limits_{j=1}^k{\Upsilon}_{j,n}^*A_j{\Upsilon}_{j,n}
&&
\sum\limits_{j=1}^k{\Upsilon}_{j,n}^*
B_j
{\Upsilon}_{j,m}
\\~\\
\sum\limits_{j=1}^k{\Upsilon}_{j,m}^*
C_j
{\Upsilon}_{j,n}
&&
\sum\limits_{j=1}^k{\Upsilon}_{j,m}^*
D_j
{\Upsilon}_{j,m}
\\~\\
I_n&&0
\\~\\
0&&I_m\end{smallmatrix}\right)
\end{matrix}
\]
\[
\begin{matrix}
=\left(\begin{smallmatrix}
\sum\limits_{j=1}^k{\Upsilon}_{j,n}^*A_j{\Upsilon}_{j,n}
&&
\sum\limits_{j=1}^k{\Upsilon}_{j,n}^*B_j{\Upsilon}_{j,m}
\\~\\
\sum\limits_{j=1}^k{\Upsilon}_{j,m}^*C_j{\Upsilon}_{j,n}
&&
\sum\limits_{j=1}^k{\Upsilon}_{j,m}^*D_j{\Upsilon}_{j,m}
\\~\\
I_n&&0\\~\\
0&&I_m
\end{smallmatrix}\right)^*
\left(\begin{smallmatrix}~0&0&
-\sum\limits_{j=1}^k{\Upsilon}_{j,n}^*{\Upsilon}_{j,n}
&0\\~\\0&0&0&
\sum\limits_{j=1}^k{\Upsilon}_{j,m}^*{\Upsilon}_{j,m}
\\~\\
-\sum\limits_{j=1}^k{\Upsilon}_{j,n}^*{\Upsilon}_{j,n}
&0&0&0\\~\\0&\sum\limits_{j=1}^k{\Upsilon}_{j,m}^*{\Upsilon}_{j,m}
&0&0\end{smallmatrix}\right)
\left(\begin{smallmatrix}
\sum\limits_{j=1}^k{\Upsilon}_{j,n}^*A_j{\Upsilon}_{j,n}
&&
\sum\limits_{j=1}^k{\Upsilon}_{j,n}^*B_j{\Upsilon}_{j,m}
\\~\\
\sum\limits_{j=1}^k{\Upsilon}_{j,m}^*C_j{\Upsilon}_{j,n}
&&
\sum\limits_{j=1}^k{\Upsilon}_{j,m}^*D_j{\Upsilon}_{j,m}
\\~\\I_n&&0\\~\\0&&I_m\end{smallmatrix}\right)
\end{matrix}
\]
\[
\begin{matrix}
=\sum\limits_{j=1}^k
\left(\begin{smallmatrix}{\Upsilon}_{j,n}^*&&0\\~\\
0&&{\Upsilon}_{j,m}^*\end{smallmatrix}\right)
\underbrace{
\left(\begin{smallmatrix}A_j&&B_j\\~\\C_j&&D_j\\~\\I_n&&0\\~\\
0&&I_m\end{smallmatrix}\right)^*\overbrace{\left(\begin
{smallmatrix}~0&&0&&-I_n&&0\\~\\~0&&0&&~0&&I_m\\~\\
-I_n&&0&&~0&&0\\~\\~0&&I_m&&~0&&0\end{smallmatrix}\right)
}^{\hat{W}_{\alpha}}\left(\begin{smallmatrix}A_j&&B_j\\~\\
C_j&&D_j\\~\\I_n&&0\\~\\0&&I_m\end{smallmatrix}\right)}_{Q_j}
\left(\begin{smallmatrix}{\Upsilon}_{j,n}&0\\~\\0&{\Upsilon}_{j,m}
\end{smallmatrix}\right)
\end{matrix}
\]
\[
\begin{matrix}
=\sum\limits_{j=1}^k\left(\begin{smallmatrix}{\Upsilon}_{j,n}^*&&0\\~\\
0&&{\Upsilon}_{j,m}^*\end{smallmatrix}\right)
Q_j
\left(\begin{smallmatrix}{\Upsilon}_{j,n}&0\\~\\0&{\Upsilon}_{j,m}\end{smallmatrix}\right)
\in\overline{\mathbf P}_{n+m}~.
\end{matrix}
\]
Thus the case of $l=\alpha$ is established.
\vskip 0.2cm

\noindent
Since $\hat{W}_{\beta}$, $\hat{W}_{\gamma}$ and $\hat{W}_{\delta}$ are
just permutations of $\hat{W}_{\alpha}$, the respective constructions
are very similar and thus omitted, and the proof is complete.
\qed
\vskip 0.2cm

\noindent
Special cases of Proposition \ref{Pn:ConvexRealizationS} for $\mathcal{P}$
and for $\mathcal{DB}$ were shown in \cite{Lewk2020a} and \cite{Lewk2020c},
respectively.
\vskip 0.2cm

\noindent
We next illustrate how, by using the above results, one can generate from
a single system, a whole collection of them. For simplicity, we address

\begin{Ex}\label{Ex:GeneratingFunctions}
{\rm
Consider the three following \mbox{$2\times 2$-valued} $\mathcal{LP}$
rational functions along with the corresponding balanced realizations,
where $a, b\in\R$ are parameters.
\begin{equation}\label{eq:ExamplePro}
\begin{matrix}
F_1(z)=&{\scriptstyle\frac{1}{a^2z}}\left(\begin{matrix}\frac{b^2}{a^2}+1&&
z-\frac{b}{a}\\~\\-(z+\frac{b}{a})&&1\end{matrix}\right)&&&
R_{F_1}=&\left({\footnotesize\begin{array}{rr|rr}0&0&\frac{1}{a}&0\\
0&0&\frac{b}{a^2}&-\frac{1}{a}\\ \hline\frac{1}{a}&\frac{b}{a^2}&0&
\frac{1}{a^2}\\0&-\frac{1}{a}&-\frac{1}{a^2}&0\end{array}}\right)
\\~\\
F_2(z)=&{\scriptstyle\frac{1}{z^2+1}}\left(\begin{matrix}a^2z&&a(bz-a)\\~\\
a(bz+a)&&(a^2+b^2)z\end{matrix}\right)&&&R_{F_2}=&\left({\footnotesize
\begin{array}{rr|rr}0&1&a&b\\-1&0&0&-a\\ \hline a&0&0&0\\b&-a&0&0
\end{array}}\right)
\\~\\
F_3(z)=&{\scriptstyle\frac{z}{1+z^2}}\left(\begin{matrix}1&&\frac{b}{a}-z
\\~\\
\frac{b}{a}+z&&\frac{b^2}{a^2}+1\end{matrix}\right)&&&R_{F_3}=&
\left({\footnotesize\begin{array}{rr|rr}0&1&1&\frac{b}{a}\\-1&0&0&1\\
\hline 1&0&0&-1\\ \frac{b}{a}&1&1&0\end{array}}\right).
\end{matrix}
\end{equation}
Following the first part of Section \ref{Sec:ExtremeCases}, each of
these three functions satisfy,
\[
\begin{matrix}
\begin{smallmatrix}-F(z)\in\overline{\mathbf L}_{I_2}&&z\in\C_L
\\~\\
F(z)\in{i}\overline{\mathbf H}_2&&z\in{i}\R
\\~\\
F(z)\in\overline{\mathbf L}_{I_2}&&z\in\C_R\end{smallmatrix}
&&&{\rm and/or}&&&\left(\begin{smallmatrix}-I_2&&0\\~\\~0&&I_2
\end{smallmatrix}\right)R_F+{R_F}^*\left(\begin{smallmatrix}-I_2&&0
\\~\\~0&&I_2\end{smallmatrix}\right)=0_{4\times 4}~.
\end{matrix}
\]
To employ Proposition \ref{Pn:ConvexRationalFunctions} to generate
additional rational functions, let now
\mbox{${\Upsilon}_j\in\C^{2\times 2}$} be arbitrary so that~
\mbox{$\sum\limits_{j=1}^3{\Upsilon}_j^*{\Upsilon}_j=I_2$}.
Then, with $F_j(s)$ from Eq. \eqref{eq:ExamplePro}, one has that~
\mbox{$\sum\limits_{j=1}^3{\Upsilon}_j^*F_j(z){\Upsilon}_j$} ~is a
$2\times 2$-valued positive real odd rational function.
\vskip 0.2cm

\noindent
Similarly, to generate additional systems by employing Proposition
\ref{Pn:ConvexRealizationS}
let now \mbox{$\tilde{\Upsilon}_j, \tilde{\tilde{\Upsilon}}_j\in\C^{2\times 2}$}
be arbitrary so that
\mbox{$\sum\limits_{j=1}^3\tilde{\Upsilon}_j^*\tilde{\Upsilon}_j=I_2$} and
\mbox{$\sum\limits_{j=1}^3\tilde{\tilde{\Upsilon}}_j^*\tilde{\tilde{\Upsilon}}_j=I_2$}.
Then, with 
$R_{F_j}$ 
from Eq. \eqref{eq:ExamplePro}, one has that
\mbox{$\sum\limits_{j=1}^3
\left(\begin{smallmatrix}\tilde{\Upsilon}_j&&0\\~\\0&&\tilde{\tilde{\Upsilon}}_j
\end{smallmatrix}\right)^*
R_{F_j}
\left(\begin{smallmatrix}\tilde{\Upsilon}_j&&0\\~\\0&&\tilde{\tilde{\Upsilon}}_j
\end{smallmatrix}\right)$} is a \mbox{$(2+2)\times(2+2)$} realization (recall,
not necessarily minimal) of a positive real odd rational function.
\vskip 0.2cm

\noindent
Finally note that we actually started from a single system $F_1$. Indeed, $F_2$
is defined as, \mbox{$R_{F_2}=(R_{F_1})^{-1}$} (in the sense of inverting a
constant $4\times 4$ matrix). Now, \mbox{$F_3(z)=\left(F_1(z)\right)^{-1}$} (in
the sense of that the product of pair of rational functions, each of degree two,
yields a zero degree rational function, i.e. $F_3(z)F_1(z)\equiv I_2$).
}
\qed
\end{Ex}
\vskip 0.2cm

\noindent
We conclude by remarking that, following Section \ref{Sec:ExtremeCases}
and Example \ref{Ex:GeneratingFunctions}, one can formally specialize
Proposition \ref{Pn:ConvexRealizationS} to balanced realizations of
$\mathcal{LB}$ and of $\mathcal{HB}_{\eta}$ functions.


\end{document}